%% file: fmn.tex
\documentclass[12pt]{amsart}

\usepackage{comment}
\input bookman.sty
\boldmath

\newcommand{\bee}{\begin{equation}}
\newcommand{\eee}{\end{equation}}

\usepackage{helvet}

\usepackage{graphics}
\usepackage{amssymb}
\usepackage{amsxtra}
\usepackage{amsmath}
\usepackage{mathrsfs}

\usepackage[arrow, matrix]{xy}
\xyoption{frame}

\theoremstyle{plain}

\newtheorem{theo}{Theorem}[section]
\newtheorem{lemm}[theo]{Lemma}
\newtheorem{prop}[theo]{Proposition}
\newtheorem{coro}[theo]{Corollary}
\newtheorem{conj}[theo]{Conjecture}

\theoremstyle{definition}

\newtheorem{defi}[theo]{Definition}
\newtheorem{rema}[theo]{Remark}

\newtheorem{exam}[theo]{Example}

\newcommand{\lb}{\label}

\newcommand{\rrf}[1]{(\ref{#1})}

\newcommand{\arrh}[3]
{
\xymatrix{
{#1} \ar[r]^<<<<{#2}  &{#3}
}
}

\newcommand{\arlh}[3]
{
\xymatrix{
{#1} & \ar[l]_<<<<{#2}  {#3}
}
}

\newcommand{\arrr}[1]
{\arrh {}{#1}{}}

\newcommand{\arr}
{\arrr {}}

\newcommand{\arrl}[1]
{\arlh {}{#1}{}}

\newcommand{\arl}
{\arrl {}}

\newcommand{\arrto}
{\xymatrix{{} \ar@{|-{>}}[r]  & {} } }

\newcommand{\arrinto}
{\xymatrix{{} \ar@{^{(}->}[r]  & {} } }

\newcommand{\becc}{\begin{comment}}                   
\newcommand{\encc}{\end{comment}}

\input auxxx.tex

\begin{document}

\title[Novikov homology and jump loci]
{Novikov homology, jump loci and Massey products}
\author{Toshitake Kohno and  Andrei Pajitnov}
\address{Kavli IPMU (WPI), Graduate School of Mathematical Sciences, the
University of Tokyo, 
3-8-1 Komaba, Meguro-ku, Tokyo 153-8914, Japan}
\email{ kohno@ms.u-tokyo.ac.jp}
\address{Laboratoire Math\'ematiques Jean Leray 
UMR 6629,
Universit\'e de Nantes,
Facult\'e des Sciences,
2, rue de la Houssini\`ere,
44072, Nantes, Cedex}                    
\email{ andrei.pajitnov@univ-nantes.fr}
\begin{abstract}
Let $X$ be a finite CW 
complex, and $\r:\pi_1(X)\to \GL(l,\cc)$ a representation.
Any cohomology class $\a\in H^1(X,\cc)$ gives rise
to a deformation $\g_t$ of $\r$ defined by
$\g_t(g)=\r(g)\exp(t\langle\a, g\rangle)$. We show 
that the cohomology of $X$ with local coefficients $\g_{gen}$
corresponding to the generic point of the curve $\g$ 
is computable from a spectral sequence starting 
from $H^*(X,\r)$. We compute the differentials of the 
spectral sequence in terms of the Massey products. We show 
that the spectral sequence degenerates in case when $X$ is a K\"ahler
manifold and $\r$ is semi-simple.

If $\a\in H^1(X,\rr)$ one associates to the triple 
$(X,\r,\a)$ the twisted Novikov homology 
(a module over the Novikov ring). We show that the twisted 
Novikov Betti numbers equal the Betti numbers of $X$ with coefficients
in the local system $\g_{gen}$. We investigate the dependence
of these numbers on $\a$ and prove that they are constant 
in the complement to a finite number of proper vector subspaces 
in $H^1(X,\rr)$.

\end{abstract}

\maketitle

\section{Introduction}
\label{s:intro}

Let $X$ be a finite connected CW-complex; denote its 
fundamental group by $G$.
Let $\r: G\to\GL(l,\cc)$ be a representation.
Any cohomology class
$\a\in H^1(X,\cc)$ gives rise to the following deformation
of $\r$:
$$
\g_t:G\to \GL(l,\cc), \ \ \g_t(g)=e^{t\langle\a, g\rangle} \r(g).
$$
The  cohomology groups of $X$
with local coefficients $\g_t$ are isomorphic 
for all $t$ except a subset containing only isolated points.
The cohomology group $H^*(X, \g_{gen})$ corresponding to the generic point
of the curve $\g_t$ is the first main 
object of study in the present paper.
We prove that there is a spectral sequence $\EEEE^*_r$
starting from the homology of $X$ with coefficients in $\r$,
converging to $H^*(X, \g_{gen})$,
and the differentials in this spectral sequence are computable in terms 
of some special higher 
Massey products with $\a$.
The first differential in this spectral sequence is  the 
homomorphism $L_\a$
of multiplication by $\a$ in the $\r$-twisted cohomology of $X$.

This type of spectral sequences appeared in the paper of S.P. Novikov
\cite{NovikovBH}
in the de Rham setting. It was generalized to the case of 
cohomology with coefficients in a field of arbitrary characteristic in 
the paper \cite{PajitnovMassey} of the second author. See also the papers
\cite{FarberTopClosedForms}, 
\cite{FarberLSforClosedForms}
of M. Farber. Our present construction
is close to the  original ideas  of S.P. Novikov. The main technical novelty
of our present approach is the systematic use of 
{\it formal exponential deformations} (see $\S\S$ \ref{s:def_diff}, 
\ref{s:spec_twis}).
This allows to avoid the convergency issues for power series,
which occur in Novikov's idea of the proof of the first main theorem of his paper
(see \cite{NovikovBH}, page 553).

If $\r$ is the trivial representation, the differentials in the spectral
sequence above are the usual Massey products
in the ordinary cohomology with slightly reduced indeterminacy:
$d_r(x) = \langle \a, \ldots , \a, x \rangle$, \ 
 see $\S$ \ref{s:def_diff}.
Thus the spectral sequence degenerates when the space $X$ is formal,
by the classical argument of 
P. Deligne, Ph. Griffiths, J. Morgan, D. Sullivan
\cite{DGMS} which applies here as well ($\S$ \ref{s:def_diff}).
Thus we have 
\begin{equation}\label{f:degen}
\Ker L_\a\big/ \Im L_\a 
\approx 
H^*(X,\g_{gen}).
\end{equation}

In case when $\r$ is not trivial the situation is more complicated.
The spectral sequence degenerates, in particular, 
when $X$ is a K\"ahler manifold and $\r$ is a 
semi-simple representation
(\S \ref{s:form_hypform}, Prop. \ref{p:high_rank_hypform}).
The proof uses C. Simpson's theory of Higgs bundles \cite{SimpsonHiggs}.
Thus the isomorphism \rrf{f:degen} holds also in this case.
We introduce a class of {\it strongly formal spaces}
for which all the spectral sequences $\EEEE^*_r$ corresponding 
to 1-dimensional representations 
degenerate in their second term.
An example from the work of H. Kasuya \cite{KasuyaMinLef}
shows that there exist formal spaces which are not strongly formal.

In the literature there are several other constructions 
of the spectral sequences related to the cohomology with twisted coefficients, such as 
the equivariant spectral sequence, introduced and studied by S. Papadima and A. Suciu 
\cite{PapadimaSuciuEqSpSeq}. The isomorphism \rrf{f:degen} was obtained also in the recent work \cite{DimcaPapadimaNonAbJL}
of A. Dimca and S. Papadima for the case when $\r$ is the trivial representation.
One of the advantages of our method is that 
it allows an explicit computation of the higher differentials of the spectral sequence
and leads to the proof of the isomorphism \rrf{f:degen} for deformations of non-trivial representations.

If $\a$ is a real cohomology class, there is another geometric construction
related to $\r$ and $\a$, namely, the {\it twisted Novikov homology} 
introduced in the works of H. Goda and the second author 
see \cite{GodaPajitnovTwiNov}, \cite{PajitnovTwiAlThurst}. 
This construction
associates to $X$, $\r$ and $\a$ a module over a corresponding 
Novikov ring
$\wh L_{m,\a}$.
The rank and torsion numbers of this module 
are called {\it the twisted Novikov Betti numbers} 
and {\it twisted Novikov torsion numbers}; they 
provide lower bounds for numbers of zeros of any Morse form
belonging to the de Rham cohomology class $\a$
(see \S \ref{s:twi_nov_B} for details).
 These invariants detect 
the fibered knots in $S^3$ as it follows from the recent work
of S. Friedl  \cite{FriedlVidussiVanTwiAl}.
For a given space $X$ and the representation $\r$
the Novikov numbers depend on
the cohomology class  $\a\in H^1(X,\rr)$. The case of the torsion numbers 
was studied in \cite{PajitnovShrp} and  \cite{PajitnovTwiAlThurst}.
It is proved there that 
the torsion numbers are constant in the open polyhedral
cones formed by finite intersections of certain half-spaces in $\rr^n$,
where $n=\rk H_1(X,\zz)$. Similar analysis applies to
the Novikov Betti numbers, which are of main interest to us in the present work.
We prove in \S \ref{s:twi_nov_B} that these numbers do not depend on $\a$
in the complement  to a finite number of proper vector subspaces.  
In general the set of all $\a$ for which 
the Novikov Betti number 
$\wh b^\r_k(X,\a)$ 
is greater by $q$ than the generic value 
(the {\it jump loci} 
for the Novikov numbers)
is a union of a finite number of proper vector subspaces,
see \S \ref{s:twi_nov_B}, Prop. \ref{p:thins_twi}.

It is known that the Novikov homology and the homology 
with local coefficients
are related to each other. This was first observed in the paper \cite{PajitnovAn}
of the second author, see also Novikov \cite{NovikovBH}.
Similar result holds also for the twisted Novikov homology,
namely,  we prove (see \S \ref{s:twi_nov_B}, Prop. \ref{p:equal})
that for $\a\in H^1(X,\rr)$
we have
$$
\wh b_k^\r(X,\a) = \b_k(X,\g_{gen}).
$$
This implies several corollaries about both families of numerical invariants.
One corollary is that for given $\r$ and $\a$ the jump loci for
$\b_k(X,\r,\a)$ are unions
of proper vector subspaces. 
On the other hand the twisted Novikov Betti numbers 
are computable from the Massey spectral sequence.integral hyperplanes
In the case of degeneracy of this spectral sequence,
they equal the dimension of its second term.

\section{Exact couples}
\label{s:ex_c}

In this section we recall the definition of the spectral sequence of an exact couple 
(following \cite{MasseyExCoup}, \cite{HuHomotTh}) and give an 
equivalent description of the successive terms of the spectral 
sequence, which will be useful in the sequel.

Let $\CCCC=(D,E,i,j,k)$ be an exact couple, so that we have an exact  triangle
$$\xymatrix{
D  \ar[rr]^i &  & D \ar[dl]^j\\
& E \ar[ul]^k  & \\
}$$
We will usually abbreviate the notation to $\CCCC=(D,E)$ and call $D$ and $E$ the first, respectively the second
component of the exact couple.
Following W. Massey we define the derived exact couple setting
$$E'=\Ker(j k)/\Im (j k),\ \  D'=i(D)$$
and defining $j', k'$  suitably.
Iterating the process we obtain a sequence of exact couples
$\CCCC_r=(D_r, E_r)$, the initial couple being numbered as $\CCCC_1$;
this sequence is called the spectral sequence associated to the exact couple $\CCCC$.

We will need an alternative description of the groups $E_r$ and the maps $j_r, k_r$.
\begin{defi}
\label{d:spec_diff}
\been\item
For $r\geq  2$ let $Z_r$ be the subgroup of all elements $x\in E$ such that
$k(x)=i^{r-1}(y)$ for some $y\in D$. We put $Z_1=E$.
\item For $r\geq  1$ let $B_r$ be the subgroup of all elements $z\in E$, such that
$z=j(y)$ for some $y\in D$ with $i^{r-1}(y)=0$.
\enen
\end{defi}
The following properties are easy to check:
$$Z_1=E\sps Z_2=\Ker(jk)\sps Z_3 \ldots \sps Z_r \sps Z_{r+1} \ldots $$
$$ B_1=\{0\}\sbs B_2=\Im(jk)\sbs B_3\sbs\ldots \sbs B_r \sbs B_{r+1}\sbs \ldots $$
$$ B_i\sbs Z_j \qquad {\rm for ~ every ~} i,j. $$

Put 
$$\wi E_r=Z_r/B_r,\ \ \  D_r=\Im i^r.$$

Define a \ho~ $\wi k_r:\wi E_r\to D_r$ 
setting $k_r(x)=k(x)$ for every $x\in Z_r$. 
Define a \ho~ $\wi j_r: D_r\to \wi E_r$ as follows: if $x\in D_r$ and $x=i^r(y)$, then
put $\wi{j_r}(x)= [j(y)]$. 
It is easy to check that these homomorphisms are well-defined and give rise to an exact couple $\wi{\CCCC_r}= (D_r, \wi{E_r})$:
$$
\xymatrix{
D_r  \ar[rr]^i &  & D_r \ar[dl]^{\wi{j_r}}\\
& \wi{E_r} \ar[ul]^{\wi{k_r}}  & \\
}$$

The proof of the next is in a usual diagram chasing:

\begin{prop}
 The exact couples $\CCCC_r$ and $\wi{\CCCC_r}$
are isomorphic for any $r$. $\qs$
\end{prop}

\section{Formal deformations of differential algebras and their spectral sequences }
\label{s:def_diff}

Let 
$$
\AA^*=\{\AA^k\}_{k\in \nn}=\{\AA^0\arrr d \AA^1\arrr d  \ldots \}
$$
be a graded-commutative differential algebra (DGA) 
over a field $\kk$ of characteristic zero.
Let $\NN^*$ be a graded differential module (DGM) over $\AA^*$
(that is, $\NN^*$ is a graded module over $\AA^*$ endowed with a differential
which satisfies the Leibniz formula \wrt~ the pairing
$\AA^*\times \NN^*\to\NN^*$). We will use the same symbol $d$
to denote the differentials in both $\AA^*$ and $\NN^*$,
since no confusion is possible.
We denote by $\AA^*[[t]]$ the algebra of formal power series over
$\AA^*$ endowed with the differential extended from the differential of
$\AA^*$. Let $\xi\in \AA^1$ be a cocycle. Consider the $\AA^*[[t]]$-module
$\NN^*[[t]]$
and endow it with the differential
$$
D_tx=dx+t\xi x.
$$
Then 
$\NN^*[[t]]$
is a DGM over $\AA^*[[t]]$, 
and we have an exact sequence of DGMs:
\begin{equation}\label{f:ex_seq}
 0 \arr \NN^*[[t]] \arrr {t} \NN^*[[t]]\arrr {\pi}  \NN^*  \arr 0
\end{equation}
where $\pi$ is the natural projection $t\arrto 0$.
The induced long exact sequence in cohomology can 
be considered as an exact couple

\begin{equation}\label{f:def_ex_c}
\xymatrix{
H^*\big(\NN^*[[t]]\big)  \ar[rr]^t &  & H^*\big(\NN^*[[t]]\big) \ar[dl]^{\pi_*}\\
& H^*( \NN^*) \ar[ul]^{\delta}  & \\
}
\end{equation}

\begin{prop}\label{p:invar}
The spectral sequence induced by the exact couple \rrf{f:def_ex_c}
depends only on the cohomology class of $\xi$.
\end{prop}
\Prf
Let $\xi_1,\xi_2\in A^1$ be cohomologous cocycles, $\xi_1=\xi_2+df$
with $f\in A^0$. Let $D_t=d+t\xi_1, ~ D'_t=d+t\xi_2$
be the corresponding differentials.
Multiplication by $e^{tf}\in A^0[[t]]$
determines an isomorphism 
$F:\NN^*[[t]] \to \NN^*[[t]]$,
 commuting with the differentials, namely, 
$F( D_t \o)=D'_t(F(\o)).$ 
Thus the exact sequences \rrf{f:ex_seq}
corresponding to $\xi_1$ and  $\xi_2$
are isomorphic, as well as the exact couples \rrf{f:def_ex_c}
and their spectral sequences. $\qs$

\begin{defi} 
\label{d:def_spec_seq}
Put $\alpha = [\xi]$. The spectral sequence associated 
to the exact couple \rrf{f:def_ex_c}
is called {\it deformation spectral sequence} and denoted 
by 
$$\EEEE^*_r(\NN^*, \alpha)=
\Big(D^*_r(\NN^*, \alpha),\  E^*_r(\NN^*, \alpha)\Big).
$$
If the couple $(\NN^*, \alpha)$
is clear from the context,
we suppress it in the notation and write 
$\EEEE^*_r$,
respectively, 
$\ D^*_r, \ E^*_r$.
\end{defi}

Denote by 
$$L_\alpha: H^*( \NN^*)\to H^*( \NN^*)$$
the  multiplication by $\alpha$.
It is clear that the first differential in the
spectral sequence equals $L_\a$ and therefore 
$$E^*_2= \Ker L_\alpha/\Im L_\alpha.$$

We are going to  compute the higher differentials in this spectral sequence
in terms of special Massey products.
Let $a\in  H^*(\NN^*)$. 
An {\it $r$-chain starting from $a$} is a sequence of elements
$\o_1, \ldots, \o_r\in \NN^*$
such that 
$$
 d\o_1=0,\ \  [\o_1]=a,\ \   d\o_2=\xi\o_1,\ \ldots,\   d\o_r=\xi\o_{r-1}.
$$
Denote by 
$MZ^m_{(r)}$
 the subspace of all 
$a\in   H^m(\NN^*)$
 such that there exists an $r$-chain
starting from $a$. Thus 
$$MZ^m_{(1)}=H^m(\NN^*), 
\ \ MZ^m_{(2)}=
\Ker \Big( L_\alpha:H^m(\NN^*)\to H^{m+1}(\NN^*) \Big). $$
Denote by 
$MB^m_{(r)}$
 the subspace of all 
$\beta\in   H^m(\NN^*)$ such that there exists an $(r-1)$-chain
$(\o_1, \ldots , \o_{r-1})$ with $\xi\o_{r-1}$ belonging to $\beta$. 
By definition 
$$MB^m_{(1)}=0,
\ \ 
MB^m_{(2)}=
\Im \Big( L_\alpha:H^{m-1}(\NN^*)\to H^{m}(\NN^*) \Big). $$
It is clear that 
$MB^m_{(i)}
\sbs
MZ^m_{(j)}$
for every $i, j$.
Put
$$
MH^m_{(r)}
=
MZ^m_{(r)} \Big/
MB^m_{(r)}.
$$
In the next definition we  omit the upper indices and write
 $MH_{(r)},MZ_{(r)} $ etc.  in order to simplify the notation.
\begin{defi}
 Let $a\in  H^*(\NN^*)$, and $r\geq 1$.
We say that the $(r+1)$-tuple Massey product
$\langle \xi, \ldots, \xi,a\rangle$
is defined, if $a\in MZ_{(r)}$.
In this case choose any $r$-chain 
$(\o_1, \ldots , \o_{r})$ 
starting from $a$. The cohomology class of $\xi\o_{r}$
is in $MZ_{(r)}$
(actually it is in $MZ_{(N)}$ for every $N$)
and it is not difficult to show that it is well defined 
modulo $MB_{(r)}$.
The image of $\xi\o_r$ in $MZ_{(r)}/MB_{(r)} $  is called the 
$(r+1)$-tuple Massey product of $\xi$ and $a$:
$$
\langle \xi, a \rangle_{(r+1)} = 
\Big\langle \ \underset{r}{\underbrace {\xi, \ldots, \xi}},\ a\ \Big\rangle
\in MZ_{(r)}\Big/MB_{(r)}.
$$
\end{defi}

\begin{exam}
 The double Massey product $\langle \xi, a \rangle_{(2)}$
equals $\xi a$, the triple Massey product $\langle \xi, a \rangle_{(3)}$
equals the cohomology class of $\xi\o_2$ where  $ d \o_2=\xi\o_1$, and $\ 
[\o_1]=a$, etc.
\end{exam}
The correspondence $a\arrto \langle \xi, a \rangle_{(r+1)}$ 
gives rise to a well-defined homomorphism of degree $1$
$$
\Delta_r: MH_{(r)} \arr  MH_{(r)}.
$$
The next proposition is proved by an easy diagram chasing argument.
\bepr\label{p:diffs_in_F}
For any $r$ we have $\Delta_r^2=0$, and 
 the  cohomology group $H^*(MH^*_{(r)}, \Delta_r)$ is isomorphic to 
 $MH^*_{(r+1)}$. $\qs$
 \enpr

\beth\label{t:compar_spec_seq}
For any $r\geq 2$ there is an isomorphism
$$
\phi:
MH^*_{(r)}\arrr \approx E^*_{r}
$$
commuting with differentials. 
\enth
\Prf
Recall from Section \ref{s:ex_c}
a spectral sequence $\wi\EEEE^*_r$
isomorphic to $\EEEE^*_r$. It is formed
by exact couples
$$
\xymatrix{
D_r  \ar[rr]^i &  & D_r \ar[dl]\\
& \wi{E_r} \ar[ul]^{\wi{\delta_r}}  & \\
}
$$
where 
$D_1=H^*( \NN^*), \ D_r=\Im\Big(t^{r-1}: D \to D\Big)$, and
$\wi E_r=Z_r/B_r$ (the modules  $Z_r, \ B_r$ are described 
in the definition 
\ref{d:spec_diff}).

\bele\label{l:mz_z}
1) $MZ_{(r)}=Z_r$, \ 2) $MB_{(r)}=B_r$.
\enle
\Prf
We will prove 1), the proof of 2) is similar.
Let $\zeta\in H^*( \NN^*)$ and  $z$ be a cocycle belonging to $\zeta$.
Then $\delta(\z)$ equals the cohomology class of $\xi z\in \NN^*[[t]]$; and 
$\zeta\in Z_r$ if and only if there is $\mu\in \NN^*[[t]]$ such that 
\begin{equation}\label{f:series}
\xi z - D_t \mu \in t^{r-1} \NN^*[[t]].
\end{equation}
This condition is clearly equivalent to the existence of 
a sequence of elements $\mu_0, \mu_1, \ldots \in \NN^*[[t]]$
such that
\begin{equation}\label{f:series2}
\xi z = d\mu_0,\ \  {\rm and }\ \  d\mu_i+\xi \mu_{i-1}=0\ \  {\rm for\ every  }\ \  0\leq i\leq r-2.
\end{equation}
The condition \rrf{f:series2} is in turn equivalent to
the existence of an $r$-chain starting from $\z$. $\qs$

The Lemma implies that  $\wi E_r^*\approx MH^*_{(r)}$
and it is not difficult to prove that this isomorphism is compatible
with the boundary operators. $\qs$

In view of the Proposition \ref{p:invar} we obtain the next Corollary.
\beco\lb{c:de_diff_invar}
Let $\xi\in \AA^1$ be a cocycle. The graded groups $MH^*_{(r)}$
defined above depend only on the cohomology class of
$\xi$, which is denoted by $\a$. $\qs$
\enco

Therefore the differentials in the spectral sequence $\{\EEEE^*_r\}$ are equal
to the higher Massey products with the cohomology class of $\xi$. 
Observe that these Massey products, defined above,  have smaller
indeterminacy than the usual Massey products.
The second term of the spectral sequence is described therefore
in terms of multiplication by the cohomology class $\a=[\xi]$. It is convenient 
to give a general definition.
\bede\lb{d:cocyle_cohom}
Let $\KK^*$ be a differential graded algebra, and $\theta$ be an element of odd 
degree $s$. Denote by $L_\t: \KK^*\to \KK^{*+s}$ the homomorphism of multiplication by
$\t$. The quotient 
$$
\Ker L_\t\Big/ \Im L_\t
$$
is a graded module which 
is denoted by $\HH^*(\KK^*, \t)$ and is called
{\it $\t$-cohomology of $\KK^*$}.
The dimension  of $\HH^k(\KK^*, \t)$ 
is denoted by $\BB_k(\KK^*, \t)$. 
\end{defi}
We have therefore
$$
E^*_2(\NN^*, \a) 
\approx
\HH^*(H^*(\NN^*), \a).
$$

Let us consider some examples, which will be important for the sequel.
\begin{exam}\label{e:algebra}
 Let $\NN^*=\AA^*$.
 The \ho~ $L_\a$ 
 is the multiplication by $\a$ in the cohomology
 $H^*(\AA^*)$, and the differentials in the spectral sequence $\EEEE_r$
 are the higher Massey products induced by the ring 
 structure in $H^*(\AA^*)$.
 \end{exam}
\begin{exam}\label{e:algebra_eta}
 Let $\AA^*$ be a DGA and $\eta\in \AA^1$ be a cocycle.
 Endow the algebra $\AA^*$ with the differential $d_\eta$
 defined by the following formula:
 $$d_\eta(x)=dx+\eta x;$$
 we obtain a DGM over $\AA$, which we will denote by
 $\wi \AA_\eta^*$. For an element $\a\in H^1(\AA^*)$
 we obtain a spectral sequence 
 $\EEEE_r^*$ with
 $$
 E_1^*= H^*(\wi\AA_\eta^*); \ \ 
  E_2^*= 
\HH^*(H^*(\wi\AA_\eta^*), \alpha).
  $$
   \end{exam}
   \begin{exam}\label{e:loc_coef}
   Let $M$ be a connected $\smo$  manifold, and $E$ be an $l$-dimensional 
   complex flat bundle over $M$.
   Denote by $\rho:\pi_1(M)\to \GL(l, \cc)$ the monodromy of $E$.
   Let $A^*(M)$ be the algebra of complex differential forms on $M$.
   The space $A^*(M,E)$ of the differential forms with 
   coefficients in $E$ is a DGM over $A^*(M)$;
   its cohomology is isomorphic to the cohomology $H^*(M,\rho)$  with 
   local coefficients \wrt~ the representation $\rho$.
For a de Rham cohomology class $\a\in H^1(M)$ 
we obtain therefore a spectral sequence 
   $\EEEE_r$ with
   $$
 E_1^*= H^*(M, \rho); \ \ 
  E_2^*= 
\HH^*(H^*(M, \rho), \alpha).
  $$

\end{exam}

Now let us consider some cases when the spectral sequences constructed above, 
degenerate in its second term.
Recall that a differential graded algebra $\AA^*$ 
is called {\it formal} if it has the same minimal model as 
its cohomology algebra. Here is a useful characterization of minimal formal
algebras. 

\beth\lb{t:form_dgms} (P. Deligne, Ph. Griffiths, J. Morgan, D. Sullivan,
\cite{DGMS}, Th. 4.1)
Let $\AA^*$ be a minimal algebra over a field of characteristic zero,
generated (as a free graded-commutative algebra) 
in degree $k$ by a vector space $V_k$; denote by $C_k\sbs V_k$
the subspace of the closed generators.
The algebra $\AA^*$ is formal if and only if 
in each $V_k$ there is a direct complement $N_k$ to $C_k$ in $V_k$,
such that any cocycle in the ideal 
generated by $\oplus_k( N_k)$ is cohomologous to zero.
\enth
This theorem leads to the proof of the well-known property that 
in formal algebras all Massey products vanish (see \cite{DGMS}). We will show that
a similar result holds for the special Massey products (Example \ref{e:algebra}).

\beth\lb{t:formal_massey}
Let $\AA^*$ be a formal differential algebra, $\alpha\in H^1(\AA^*)$.
Then the spectral sequence $\EEEE^*_r(\AA^*, \alpha)$
degenerates at its second term, and  
$$E^*_2(\AA^*, \alpha)=
\HH^*(H^*(\AA^*), \a).$$
\enth
\Prf
It suffices to establish the property for the case
of formal minimal algebras. Let $\AA^*$ be a formal minimal algebra
with the space of generators $V_k$ in dimension $k$ decomposed as
 $V_k=C_k\oplus N_k$, see \ref{t:form_dgms}. We will prove that
$MZ_{(2)}= MZ_{(3)}= \ldots =MZ_{(r)}$
and $MB_{(2)}=MB_{(3)}= \ldots = MB_{(r)}$ for every $r\geq 2$. 

Choose $\xi\in A^1$ representing $\alpha\in H^1(A^*)$.
Let $a\in MZ_{(r)}$, and $(\o_1, \ldots, \o_r)$ be an $r$-chain starting 
from $a$. so that
$$
d\o_1=0,\ \  [\o_1]=a,\ \  d\o_2=\xi\o_1,\ \ldots, d\o_r=\xi\o_{r-1}.
$$
Denote by $\L(C_*)$ the algebra generated by the space of closed generators,
so that 
$\MM=\L(C_*)\oplus I(N_*)$. Write $\o_r=\o_r^0+\o^1_r$ with
$\o_r^0\in \L(C_*),\ \  \o^1_r \in I(N_*)$.
Then $\xi\o_{r-1}=d\o^1_r$, and 
$$
d(\xi\o^1_{r})= d(\xi\o_{r})= \xi^2\o_{r-1}=0,$$
so that $\xi\o^1_{r} $ is a cocycle belonging to $I(N_*)$.
Therefore $\xi\o^1_{r}=d \o_{r+1}$ for some $\o_{r+1}\in I(N_*)$
and we obtain an $(r+1)$-chain starting from $a$, so that $a\in MZ_{(r+1)}$.

A similar argument shows that $MB_{(2)}=MB_{(r)}$ for every $r\geq 2$. 
Therefore the spectral sequence degenerates at its second term. $\qs$

The next proposition gives a sufficient condition for the 
degeneracy of the spectral sequence associated with a differential graded
module  over a DG-algebra 
$\AA^*$.
It will be used in Section \ref{s:form_hypform}
while studying the case of K\"ahler manifolds.

\bede\label{d:formal_module}
A DG-module $\NN^*$ over a DG algebra $\AA_*$ will be  called {\it formal}
if it is a direct summand of a formal DGA $\BB^*$ over $\AA_*$,
that is, 
\begin{equation}\label{f:summa}
\BB^*=\NN^*\oplus\KK^*,
\end{equation}
where both $\NN^*$ and $\KK^*$ are differential graded $\AA_*$-submodules of $\BB_*$.
\end{defi}

\bepr\label{p:degen_module}
Let $\NN_*$ be a formal DG-module over $\AA_*$,
and $\a\in H^1(\AA^*)$.
Then the spectral sequence
$\EEEE_r^*(\NN^*, \a)$
degenerates at its second term.
\enpr
\Prf
The direct sum decomposition \rrf{f:summa} implies that
$$
\EEEE_r^*(\BB^*, \a)
=
\EEEE_r^*(\NN^*, \a)\oplus\EEEE_r^*(\KK^*, \a);
$$
the spectral sequence 
$\EEEE_r^*(\BB^*, \a)$
degenerates by the previous proposition and the result follows. 
$\qs$

\section{A spectral sequence converging to the twisted cohomology}
\label{s:spec_twis}

Let $X$ be a finite CW-complex, put $G=\pi_1(X)$. 
We endow the universal covering $\wi X$ with the natural left action of $G$,
so that  the cellular chain complex of $\wi X$ is a free finitely generated chain complex 
over $\zz G$.
Let $B$ be an integral domain and $\r$ be a left action of $G$ on the free $B$-module
$B^l$ (or, equivalently, a representation
$\rho:G\to  \GL(l, B)$).
The cohomology of the cochain complex
$$C^*(X,\rho)=\Hom_G\Big(C_*(\wi X), B^l\Big)$$
%(where $\wi X$ stands for the universal covering of $X$)
is a $B$-module 
called 
{\it the twisted cohomology of $X$ \wrt~ $\r$} and denoted by
$H^*(X,\r)$. Denote by $\{B\}$ the fraction field of $B$.
The dimension over $\{B\}$ of the localization $H^k(X,\r)\otimes \{B\}$
will be called the {\it $k$-th cohomological Betti number of $X$
\wrt~ $\r$ } and denoted by $\b^k(X,\r)$.

Let us start with a given  representation $\r:G\to  \GL(l, \cc)$.
Pick a cohomology class 
 $\a\in H^1(X,\cc)$ and consider the exponential
deformation of $\r$:
\begin{equation}\label{f:exp_def}
\g_t:G\to \GL(l, \cc), \ \  
\g_t(g)=\r(g)e^{t\langle\a, g\rangle} \ \ \ (t\in\cc).
\end{equation}
Denote by $\HH$ the ring of all entire holomorphic functions on $\cc$
and let $\L=\cc[[t]]$; we have a natural inclusion $i:\HH\arrinto \L$.
The formula \rrf{f:exp_def} 
defines a family of representations of $G$:

\been\item 
For a fixed $t\in\cc$ a  representation $\g_t:G\to \GL(l,\cc)$.
\pa\item
a  representation $\bar\g \ : \ G\to \GL(l, \HH)$
(the {\it holomorphic exponential deformation} of $\r$).
\pa\item
a  representation $\wh\g = i\circ \bar\g \ : \ G\to \GL(l, \L)$
(the {\it formal exponential deformation} of $\r$).
\enen
The inclusion $i$ extends to the inclusion of the  fields of fractions
$\{\HH\}\arrinto \{\L\}$, therefore 
$$\b^k(X,\bar\g)
=
\b^k(X,\wh\g)
$$
for every $k$.

\bele\label{l:generic_t}
For every $k$ and $t$ we have 
$\b^k(X,\g_t)\geq \b^k(X,\bar\g)$.
There is a subset $S\sbs \cc$ consisting of isolated points,
such that 
\begin{equation}\label{f:generic_tt}
\b^k(X,\g_t)
=
\b^k(X,\bar\g)
\end{equation}
for every $k$ and every $t\in \cc\sm S$.
\enle
\Prf
Let $n_k$ be the number of $k$-cells in $X$.
The boundary operator $\partial_k$ in $C_*(\wi X)$ is represented by an 
$(n_{k-1}\times n_{k})$-matrix
with coefficients in $\zz G$.
The chain complex 
$$C^*(X, \bar\g) = \Hom_G\Big(C_*(\wi X), \ \HH^l\Big)$$
computing the cohomology of $X$ with 
coefficients in the local system defined by 
$\bar\g$,  has $l\cdot n_k$ free generators in degree $k$;
its boundary operator $\delta_{k+1}:C^k\to C^{k+1}$ is given by the formula 
$$
\delta_k
=
\bar\g(\pr_{k+1}^{\rm T}),$$
(we denote by $M^T$ is the transpose of $M$).
Denote by $\r_{k+1}$ the maximal rank of non-zero minors of this matrix.
For $t\in\cc$ denote 
by $\r_{k+1}(t)$ the maximal rank of non-zero minors of the matrix 
$\bar\g\big(\pr_{k+1}^{\rm T}(t)\big).$
Then 
\begin{gather*}
\b^k(X, \bar\g)=l\cdot n_k-\r_k-\r_{k+1},\\
\b^k(X, \g_t)=l\cdot n_k-\r_k(t)-\r_{k+1}(t).
\end{gather*}
Observe that  $\r_k(t)\leq \r_k$ for every $t$, and the set of 
$t\in\cc$ where these inequalities are strict, consists of isolated points
(since the  minors of $\bar\g(\pr_{k})$ are holomorphic 
functions of the variable $t\in\cc$). The lemma follows. $\qs$

Thus
$$
\b^k(X,\wh\g)
=
\b^k(X,\bar\g)
=
\b^k(X,\g_t)
\qquad
{\rm for }\  \  t\notin S.
$$
The exact sequence
\begin{equation*}
%\label{f:ex_seq_twi}
0 \arr \L \arrr {t} \L \arr \cc  \arr 0
\end{equation*}
gives rise to a short exact sequence of complexes 
\begin{equation*}
%\label{f:ex_seq_twi}
0 \arr C^*(X,\wh\g) \arrr {t} C^*(X,\wh\g) \arr C^*(X,\rho)  \arr 0
\end{equation*}
inducing a long exact sequence in cohomology,
which can be interpreted as 
an exact couple
\begin{equation}\label{f:ex_c_twi}
\xymatrix{
H^*\big(X,\wh\g \big)  \ar[rr]^t &  & H^*\big(X,\wh\g\big) \ar[dl]\\
& H^*(X, \r) \ar[ul]  & \\
}
\end{equation}
Denote by $\WWWW^*_r(X,\r, \alpha)=(U_r^*, W_r^*)$
the induced spectral sequence. 
%It is not difficult to prove the following proposition
\bepr\label{p:sp_seq_twi}
Let $X$ be a finite connected CW-complex. Then 
$$
\dim_\cc W_\infty^k
=
\b^k(X,\wh\g)
=
\beta^k(X, \g_t) \ \ 
{\rm for ~ generic ~ }  t.
$$
\enpr
\Prf
The ring $\L$ is principal, and the $\L$-module
$H^*\big(X,\wh\g \big)$
is a finite direct sum of cyclic modules.
Only the free summands survive to $E_\infty$
and the contribution of each such summand to 
$\dim_\cc W_\infty^k$
equals $1$. $\qs$
 
The cohomology with local coefficients arises naturally in the de Rham theory.
Later on we shall work also with the {\it homology} with local coefficients, which appears
in the Morse-Novikov theory. 
We will now explain the relation between the two 
constructions.

An antihomomorphism $\t:G\to GL(l,B)$ will be also called 
{\it a right representation}; it induces a structure of a right $\zz G$-module on 
$B^l$. 
The homology of the chain complex
$$C_*(X,\t)=B^l\tens{\t} C_*(\wi X) $$
%(where $\wi X$ stands for the universal covering of $X$)
is a $B$-module 
called 
{\it the twisted homology of $X$ \wrt~ $\t$} and denoted by
$H_*(X,\t)$. 
The dimension over $\{B\}$ of the localization $H_k(X,\t)\otimes \{B\}$
will be called the {\it $k$-th homological Betti number of $X$
\wrt~ $\t$ } and denoted by $\b_k(X,\t)$.

To each representation $\r:G\to GL(l,B)$ we associate the 
right representation $\r^*:G\to GL(l,B)$, where the matrix $\r^*(g)$
is the transpose of $\r(g)$. The representation $\r^*$ will be called
{\it conjugate } to $\r$. Let $E$ be the 
free finitely generated $B$-module endowed with the structure of a left $\zz G$-module
via $\r$, then the right representation $\r^*$ corresponds to the 
right  $\zz G$-module $\Hom_B(E, B)$, which is denoted by $E^*$.
Observe that $(\r^*)^*= \r$ and $(E^*)^*\approx E$.
\bele\label{l:hom_cohom}
Let $\r:G\to GL(l,B)$ be a representation. Then
there is a natural isomorphism
\begin{equation}\label{f:hom_cohom}
\Hom_B\big(C_*(X,\rho^*), B\big)
\approx
C^*(X, \r).
\end{equation}
\enle
\Prf
Put $R=\zz G$.
For any $R$-module $B$ 
there is  a canonical isomorphism
$$
\Hom_B\big(M\tens{R}B^l, B\big)
\approx
\Hom_R\big(M, \Hom_B(B^l,B)\big).
$$
The lemma follows. $\qs$

The following Corollary is immediate.
\beco\label{c:hom_cohom}
\begin{equation}\label{f:betti_hom_cohom}
\b^k(X,\r)
=
\b_k(X,\r^*).
\end{equation}
\enco
\bere
Observe that for $l=1$ every representation $\r$ is equal to its
conjugate and therefore 
$\b^k(X,\r)
=
\b_k(X,\r)$
for 1-dimensional local systems.
\enre

\section{Comparing the two spectral sequences}
\lb{s:compar_manif}

We are interested in the formal deformations of differential algebras 
related to topological spaces, mainly  manifolds. 
Let $M$ be a connected $\smo$ manifold, denote $\pi_1(M)$ by $G$, and let
$\r:G\to \GL(l,\cc)$ be  a representation. 
Let $\wi M$ be the universal covering of $M$ and define a flat bundle $E$ on $M$
as follows:
$$
E=(M\times \cc^l)\big/ \sim\ \ 
{\rm where \ } 
(gm, \xi)\sim (m, g^{-1}\xi)
{\rm \ \ for } \ g\in G.
$$
%Then the sections of $E$ are  $G$-equivariant functions $\wi M  \to \cc^l$.
The module  $\NN^*=A^*(M, E)$  of $E$-valued 
$\smo$ differential forms on  $M$ is a 
DGM over the algebra $A^*(M)$ of $\cc$-valued $\smo$
differential forms on $M$.
The elements of $\NN^*$ can be considered as $G$-equivariant differential forms on
$\wi M$ with values in $\cc^l$.
Let $\alpha\in H^1(M,\cc)$. The corresponding exact couple

\begin{equation}\label{f:def_ex_c_manif}
\xymatrix{
H^*\Big(\NN^*[[t]]\Big)  \ar[rr]^t &  & H^*\Big(\NN^*[[t]]\Big) \ar[dl]^{\pi_*}\\
& H^*(\NN^*) \ar[ul]  & \\
}
\end{equation}
gives rise to the deformation spectral sequence
$\EEEE^*_r(\NN^*,\alpha)$
(see \S \ref{s:def_diff}).
Let $\wh\g$ be the formal exponential deformation
of $\r$ corresponding to the class $\a$:

$$
\wh\g(g)=\r(g)e^{t\langle \a, g\rangle} 
\in
\GL(l,\L), \ \ \L=\cc[[t]].
$$
We associate to this deformation the spectral 
sequence $\WWWW^*(M,\r,\a)$
(see \S \ref{s:spec_twis}).

\beth\label{t:compar_manif}
The spectral sequences 
$\EEEE^*_r(\NN^*,\alpha)$
and 
$\WWWW^*_r(M,\r,\a)$
are isomorphic.
\enth
\Prf
The $A^*(M)$-module $A^*(M,E)[[t]]$ 
can  be considered as the vector space 
of  exterior differential 
forms on $M$ with 
values in $E[[t]]$.
%Let $\pi:\wi M\to M$ be the universal covering of $M$.
%The vector bundle $E$ is isomorphic to
%$\wi M\underset{\rho}{\times} \cc^l$, and
%$\pi^*E\approx \wi M\times \cc^l$.
Denote by $\TT^*(M)$
the $A^*(M)$-submodule
of  $A^*(\wi M, \cc^l[[t]])$,
consisting of  differential forms on $\wi M$
which are equivariant with respect to the representation $\wh\g$.
Choose a closed 1-form $\xi$ within the cohomology class $\a$. 
Let $F:\wi M \to\cc$ be a $\smo$ function such that
$\pi^*\xi=dF$.
The next lemma is obvious.
\bele\label{l:equiv_versus_deform}
The homomorphism
$$
\Phi: A^*(M,E)[[t]] \arr \TT^*(M)
$$ 
defined   by \ \ 
$
\Phi(\o)= e^{tF}\pi^*(\o)
$
\ \ 
is an isomorphism of DG-modules over $A^*(M)$. $\qs$
\enle

For a manifold $N$ we 
denote by $\SS_*(N)$ the 
graded group of  singular $\smo$-chains on $N$.
The integration map determines a \ho~
$$
I:
\TT^*(M)
\arr
\Hom_{G}\Big(\SS_*(\wi M), \cc^l[[t]]\Big).
$$
Here $\cc^l[[t]]$ is endowed with the structure of a $G$-module
via the representation $\wh \g: G\to \GL(l,\L)$.
\bepr\label{p:derham}
The induced map in the cohomology groups
$$
I_*:
H^*\big(\TT^*(M)\big)
\arr
H^*(M, \wh \g)
$$
is an isomorphism.
\enpr
\Prf
The argument follows the usual sheaf-theoretic proof of the 
de Rham theorem. Consider the sheaf $\TT^k$ on $M$
whose sections over an open subset $U\sbs M$
 are $\wh \g$-equivariant
$k$-forms on $\pi^{-1}(U)$ with values in $\cc^l[[t]]$. 
Denote by 
$\AA$ the sheaf on $M$, whose sections 
over $U$ are $\wh \g$-equivariant 
locally constant functions 
$\pi^{-1}(U)\to \cc^l[[t]]$.
We have an exact sequence of sheaves
\begin{equation}\lb{f:r_one}
\RR_1=\{ 0 \arr \AA\arr \TT^0\arrr d  \TT^1 \arr \ldots \}.
\end{equation}
Consider the sheaf $\ZZ^*$ of 
$\wh \g$-equivariant singular cochains:
$$
\ZZ^*(U)
=
\Hom_{G}\Big(\SS_*(\pi^{-1} (U)), \cc^l[[t]]\Big).
$$
The sequence of sheaves
\begin{equation}\lb{f:r_two}
\RR_2=\{0 \arr \AA\arr \ZZ^0\arrr \delta  \ZZ^1 \arrr\delta \ldots \}
\end{equation}
(where $\delta$ is the coboundary operator on singular cochains)
is also exact, and we have  two soft 
acyclic resolutions of the sheaf $\AA$.
The integration map $I$ induces a homomorphism $\RR_1\to \RR_2$
which equals identity on $\AA$. 
The standard sheaf-theoretic result implies that
the induced \ho~ in the cohomology groups of the complexes 
of global sections
is an isomorphism (see for example \cite{WellsAC}, Corollary 3.14). 
%The proof of the Proposition is complete. 
$\qs$

The isomorphism 
$$
I_*\circ \Phi_*: 
H^*\Big(A^*(M,E)[[t]]\Big)
\to
H^*(M, \wh \g)
$$
induces an isomorphism 
of the exact couple

\begin{equation}\label{f:compar_ex_c}
\xymatrix{
H^*\Big(A^*(M,E)[[t]]\Big)  \ar[rr]^t &  & H^*\Big(A^*(M,E)[[t]]\Big) \ar[dl]^{\pi_*}\\
& H^*(M, E) \ar[ul]  & \\
}
\end{equation}
to the exact couple
\begin{equation}\label{f:compar_ex_cc}
\xymatrix{
H^*\big(M,\wh\g \big)  \ar[rr]^t &  & H^*\big(M,\wh\g\big) \ar[dl]^{\pi_*}\\
& H^*(M, E) \ar[ul]  & \\
}
\end{equation}
and therefore the isomorphism of the spectral sequences
$$
\EEEE^*_r(\NN^*, \alpha)
\approx
\WWWW^*_r(M,\r, \alpha). \hspace{2cm} 
$$
The proof of Theorem 
\ref{t:compar_manif}
is complete.
$\qs$

Now let us proceed to general topological spaces.
The rational homotopy theory of D. Sullivan 
(see \cite{SullivanInfComp},
\cite{DGMS},
\cite{BG})
associates to 
a connected topological space 
$X$ a minimal algebra $\MM^*(X)$ over $\cc$,
well defined up to isomorphism. 

Let $\alpha\in  H^1(X,\cc)$; 
we obtain a spectral sequence 
$\EEEE^*_r(\MM^*(X), \alpha)$. It is not 
difficult to see that for the case when $X$ is a $\smo$ manifold,
this spectral sequence is isomorphic 
to the one considered in the previous section.
Indeed, we have a homotopy equivalence
$\phi : \MM^*(X)\arr A^*(X)$; the induced homomorphism of the spectral sequences
$$\EEEE^*_r(\MM^*(X),\alpha)\to 
\EEEE^*_r(A^*(X),  \phi_*(\alpha) )$$
is an isomorphism in the second term $\EEEE_2$.
Thus the two spectral sequences are isomorphic.

\beth\label{t:compar_cw}
Let $X$ be a finite CW complex. 
Then
$$
\EEEE^*_r(\MM^*(X),\alpha)
\approx
\WWWW_r^*(X, \alpha).
$$
\enth
\Prf
Let $f:X\to M$ be a homotopy equivalence of $X$ to a $\smo$ manifold $M$ 
(possibly non-compact). 
Denote by $\alpha'$ the $(f^{-1})^*$-image of $\alpha$.
The map $f$ induces a homotopy equivalence
$F: \MM^*(X) \to A^*(M)$ of DGAs.
We obtain isomorphisms of 
the spectral sequences:
\begin{gather*}
\EEEE^*_r(\MM^*(X), \alpha) \approx \EEEE^*_r\Big(A^*(M),  (f^{-1})^*(\alpha) \Big)\\
\WWWW_r^*(X,  \alpha) \approx \WWWW_r^*(M),  \alpha).
\end{gather*}
Now apply Theorem 
\ref{t:compar_manif} and the result
follows. $\qs$

\bere\label{r:simpli_derham}
Theorem \ref{t:compar_cw}
can be generalized to the case of local coefficients
so as to obtain a result similar to Theorem \ref{t:compar_manif}
in the setting of general topological spaces.
\enre

\section{Strongly formal manifolds}
\label{s:form_hypform}

The next theorem follows immediately from Theorems \ref{t:compar_cw}
and \ref{t:formal_massey}.

\beth\label{t:formal_top}
Let $X$ be a finite connected CW-complex.
Assume that $X$ is formal. Then the spectral sequences
$$
\EEEE^*_r(\MM^*(X),\a)\approx \WWWW^*(X,\a)
$$
degenerate at their second term.
\enth
Thus the dimension of the homology of $X$ with coefficients 
in a generic point of the exponential deformation of the trivial 
representation can be computed from the multiplicative
structure of the ordinary homology.
Namely, let $\a\in H^1(X,\cc)$ and let $\g_t$ be the exponential deformation
of the trivial representation:
$$
\g_t(g)=e^{t\langle \a, g\rangle }.
$$
Denote by $\b^k(X,\g_{gen})$ the Betti number of $X$ with coefficients in a
generic point of the curve $\g_t$, and let $L_\a$ be the operation of multiplication
by $\a$ in $H^k(X,\cc)$. The theorem above implies that
$$
\b^k(X,\g_{gen})=\dim_\cc
\HH^k(H^*(X,\cc), \a).
$$
The case of the spectral sequence $\WWWW^*(X,\r,\a)$
where $\r:\pi_1(X)\to GL(l,\cc)$
is a non-trivial representation, is more complicated and to guarantee
the degeneracy of the spectral sequence a stronger condition is necessary.
Let us first consider the case of 1-dimensional representations $\r$.
Let $M$ be a connected $\smo$ manifold.
Denote by $G$ the fundamental group of $M$,
let $Ch(G)$ be the group of \ho s  $G\to \cc^*=\GL(1,\cc)$.
For a character $\r\in Ch(G)$ 
denote by $E_\r$ the corresponding flat vector bundle over $M$.
Put
$$
\bar A^*(M) = \underset{\r\in Ch(G)}\bigoplus A^*(M,E_\r).
$$
The pairing $E_{\r}\otimes E_{\eta}\approx E_{\r\eta}$
induces a natural structure of a differential graded algebra
on the vector space $\bar A^*(M)$.
\begin{defi}\label{d:hyp_form}
A $\smo$ manifold $M$ is {\it strongly formal} if the algebra $\bar A^*(M)$
is formal.
\end{defi}
\beth\label{t:hyp_formal_degen}
Let $M$ be a 
strongly formal manifold, $\r\in Ch(G)$ and $\a\in H^1(M,\cc)$.
Then the spectral sequences
$$
\EEEE^*_r(\MM^*(M),\r,\a)\approx \WWWW^*_r(M,\r,\a)
$$
degenerate in their second term.
\enth
\Prf
The DG-module $A^*(M, E_\r)$ is formal;
apply Proposition \ref{p:degen_module}
and the proof is over. $\qs$

Denote by $b_k(M,\r,\g_{gen})$ the $k$-th Betti number 
of $M$ with coefficients in a generic
point of the curve 
$$
\g_t(g)=\r(g)
e^{t\langle \a, g\rangle }.
$$
\beco\lb{c:betti_hf}
Let $M$ be a strongly formal  manifold, $\r\in Ch(G)$, and $\a\in H^1(M,\cc)$.
Then 
$$
b_k(M,\r,\g_{gen}) = \BB_k(H^*(M,\r), \a).
$$
\enco
An big class of  examples of strongly formal  spaces is formed
by K\"ahler manifolds, as it follows from C. Simpson theory 
of Higgs bundles \cite{SimpsonHiggs}.
\beth\lb{t:kaeh_hypform}
Any compact K\"ahler manifold is strongly formal.
\enth
\Prf
Let $\r\in Ch(\pi_1(M))$.
The flat bundle $E_\r$ has a unique structure of a harmonic Higgs bundle
(see \cite{ArapuraHGL}, \cite{SimpsonHiggs});
the exterior differential $D_\r$ in the DG-module 
$A^*(M, E_\r)$
writes therefore as
$D_\r=D'_\r+D''_\r$, and the natural \ho s
of DG-modules 
$$
\xymatrix{
&\Big(\Ker(D'_\r), D''_\r\Big) \ar[dl] \ar[dr] & \\
\Big(A^*(M,E_\r), D_\r\Big)  &  & \Big(H^*_{DeRham}(M, E_\r), 0\Big)
}
$$
induce isomorphisms in cohomology 
(see \cite{SimpsonHiggs}, Lemma 2.2 (Formality)).
Denote the DG-module 
$(\Ker(D'_\r), D''_\r)$
by $K_\r(M)$,
and put
$$
\KK^*(M)
=
\underset{\r\in Ch(G)}\bigoplus K_\r(M).
$$The multiplicativity properties of Higgs bundles
imply that $\KK^*$ is a DG-algebra
and we have the maps of DGAs:
$$
\xymatrix{
& \KK^*(M)\ar[dl] \ar[dr] & \\
\bar A^*(M)  &  & 
\underset{\r\in Ch(G)}\bigoplus H^*(M, E_\r)
}
$$
both inducing isomorphisms in cohomology. 
The theorem follows. $\qs$

\bere\label{r:exampleNotStrForm}
There are manifolds which are formal but not
strongly formal. The example described below 
was indicated to the authors by H. Kasuya.

H. Sawai \cite{SawaiLieLatt}
constructed an 8-dimensional solvmanifold, having several remarkable properties.
H. Sawai's construction starts with a 7-dimensional solvable Lie algebra
$\mathfrak g$ generated by \break 
$A, X_1, X_2, X_3, Z_1, Z_2, Z_3$, with
\begin{gather*}
[X_1, X_2]= Z_3,\ \  [X_2, X_3]= Z_1,\ \  [X_3, X_1]= Z_2, \\
[A, X_1]=-a_1X_1,\ \  [A, X_2]=-a_2X_2, \ \ [A, X_3]=-a_3X_3,  \\  
[A, Z_1]=a_1Z_1,\ \  [A, Z_2]=a_2Z_2, \ \ [A, Z_3]=a_3Z_3,  \\
\end{gather*}
where $a_{1}, a_{2}, a_{3}$ are distinct real numbers.
(This is a generalization of the Lie algebra constructed by 
Benson and Gordon in \cite{BensonGordonKS}.) Denote by $G$ the corresponding simply connected Lie group.
H. Sawai proves that for some choice of $a_1, a_2, a_3$
there is a lattice $\Gamma$  in $G\times \RRR$, and the quotient 
is a formal space, which has a symplectic structure and 
satisfies the hard Lefschetz property,
but admits no K\"ahler structure.

The cohomology of $(G/\Gamma)\times S^1 $ with local coefficients 
in one-dimensional local systems was studied by H. Kasuya in \cite{KasuyaMinLef}.
By the Mostow theorem \cite{MostowCohom} 
the computations can be carried out in the cohomology 
of the Lie algebra $\Gamma\times\RRR$ with coefficients in 1-dimensional modules.
H. Kasuya gives an example of non-vanishing 
triple Massey product in the homology of $\Gamma\times\RRR$
with twisted coefficients. Thus $M$ is 
a formal but not a strongly formal space.

H. Kasuya informed the authors that he has 
constructed a 4-dimensional solvmanifold
which is formal but not strongly formal.
\pa

\begin{conj}\label{conj:nonstronglyformal}
For every $n\geq 2$ there exists a solvmanifold $M$, a character $\r:\pi_1(M)\to \cc^*$ and 
a cohomology class $\alpha\in H^1(M, \rr)$
such that $M$ is formal, but the differential $d_n$ in the spectral sequence 
$\WWWW_r^*(M, \r, \a)$ is non-zero.
\end{conj}

Consider now the case of representations of higher rank.
\bepr\lb{p:high_rank_hypform}
Let $M$ be a compact K\"ahler manifold
and $\r:\pi_1(M)\to \GL(l,\cc)$ be a semi-simple representation.
Then the differential graded $A^*(M)$-module $A^*(M, E_\r)$ is formal.
\enpr
\Prf 
By \cite{SimpsonHiggs}, Theorem 1 there is a harmonic metric
on the bundle $E_\r$. 
The tensor powers of this metric provide harmonic metrics on the 
bundles $E_\r^{\otimes n}$ for any $n\geq 1$.
Put
$$
\LL_\r^*(M)= \bigoplus_{n=0}^\infty A^*(M,E_\r^{\otimes n})
$$
(where  $A^*(M,E_\r^{\otimes 0}) = A^*(M)$ by convention).
Then $\LL_\r^*(M)$
is a DG-algebra. The same argument as in the proof of 
Theorem \ref{t:kaeh_hypform}
implies that this algebra is formal, and it remains to
observe that $A^*(M, E_\r)$ is a direct summand of 
$\LL_\r^*(M)$. $\qs$

\end{rema}

\section{Chain complexes over 
Laurent polynomial rings and their localizations}
\label{s:betti_nov}

In this section we discuss the Betti numbers of complexes 
over Laurent polynomial rings in view of further applications to
the Novikov Betti numbers and the homology with local 
coefficients.

Let $T$ be an integral domain and 
$
\{ T \}
$
be its fraction field. Let 
$
C_*$
be a finite free chain  complex over $T$:
$$
C_*
=
\{ 0 \arl C_0 \arl \ldots \arrl {\pr_k} 
C_k \arrl {\partial_{k+1}}  C_{k+1} \arl \ldots \}
$$
The tensor product of this  chain complex
with $\{T\}$
will be denoted by 
$\ove C_*$
and the boundary operator in 
$\ove C_*$ 
will be denoted
by $\ove \partial_*$.
The Betti number $b_k(\ove C_*)$ 
is equal to 
$$
\dim_{\{T\}} \ove C_k - \rk \ove\pr_k - \rk \ove\pr_{k+1}
=
\rk  C_k - \rk \pr_k - \rk \pr_{k+1}
$$
(where 
$
\rk \pr_k
$
stands for the maximal rank of a non-zero
minor of the matrix of $\pr_k$).

Let $\phi: T\to U$ be a \ho~ to another integral domain.
Denote by $b_k(C_*, \phi)$ the $k$-th Betti number of 
the tensor product 
$
C_*\otimes_\phi \{U\}$.
Then 
$
b_k(C_*, \phi)
\geq
b_k(C_*).
$
The inequality is strict 
if the $\phi$-images of all the 
$\rk \pr_k$-minors of 
$\pr_k$ vanish, or
if the 
$\phi$-images of all the 
$\rk \pr_{k+1}$-minors of 
$\pr_{k+1}$ vanish.
In general, for $q>0$ the condition
$
b_k(C_*, \phi)
\geq
b_k(C_*)+q
$
is equivalent to the existence of a 
number $i$, with $0\leq i \leq q$
such that all the
$(\rk \pr_k - i )$-minors of 
$\phi(\pr_k)$ 
vanish 
and all the 
$(\rk \pr_{k+1} -  (q-i))$-minors of
$\phi(\pr_{k+1})$ 
vanish.

We are interested in the case 
of Laurent polynomial rings.
Let $R$ be an integral domain and $C_*$
be a finite free chain complex over 
$
L_n=R[t_1^\pm, \ldots , t^\pm_n] = R[\zz^n].
$
Let
$
p:\zz^n\to\zz^m 
%\ \  \l:\zz^n\to R^\bu
$
be a group homomorphism. Extend it to a ring \ho~
$L_n\to L_m$
which will be denoted by the same letter $p$.
Denote by 
$\QQ_m$ is the field of fractions of $L_m$, that is, 
 the field of the rational functions 
in $m$ variables with  coefficients
in the fraction field of $R$.
Form the chain complex
$
C_*\otimes_p \QQ_m,
$
and denote by $b_k(C_*, p)$ the dimension of
the vector space 
$
H_k(C_*\otimes_p\QQ_m)
$
over $\QQ_m$.
Observe that if $p$ is injective, then 
$b_k(C_*)=b_k(C_*,p)$.
We will now study the dependance of 
$b_k(C_*, p)$
on $p$.

\bede\label{d:full}
A subgroup $G\sbs \zz^n$ is called {\it full} 
if it is a direct summand of 
$\zz^n$.
We say that a \ho~ 
$p:
\zz^n
\to
\zz^m
$
 is  {\it subordinate } to a full subgroup
 $G\sbs \Hom(\zz^n,\zz)$ and we write $p \sqsubset G$,
if all the coordinates of $p$ 
are in $G$.
\end{defi}
\bere\label{r:remk}
Let $G$ be a full subgroup of $\Hom(\zz^n,\zz)$.
Denote by $K$ the subgroup of $\zz^n$
dual to $G$. Then $p\sqsubset G$
if and only if $p~|~K = 0$.
\end{rema}

\beth\label{t:thins}
Let $C_*$ be a finite free complex over 
$L_n$. 
Let $k\geq 0, q>0$.
Then there is a finite family of proper 
full subgroups $G_i\sbs \Hom(\zz^n, \zz)$
such that for $p\in\Hom(\zz^n, \zz^m)$
the condition
$$b_k(C_*,p)\geq b_k(C_*)+q
$$
is equivalent to the following condition:
\
$p\sqsubset G_i$ for some $i$.
\enth
\Prf
Let us do the case $q=1$, the general case is similar.
Let
$\EEEE$ denote the set of all the
$(\rk\pr_k)$-minors of the matrix 
$\pr_k:C_k\to C_{k-1}$, and  all
 the
$(\rk\pr_{k+1})$-minors of the matrix 
$\pr_{k+1}:C_{k+1}k\to C_{k}$.
Let $\D\in \EEEE$, write
$\D=\sum_{g\in\zz^n} r_g \cdot g$
(where $r_g\in R$).

According to our previous observation
it suffices to study
the set $\Sigma$ of all homomorphisms
$p:\zz^n\to\zz^m$
such that $p(\D)=0$.
Let $\G=\supp\D$,
which is a finite subset of $\zz^n$.
Any \ho~ $p:\zz^n\to\zz^m$
with $p(\D)=0$
must be non-injective on $\G$.
To describe the set of all such \ho s
let us say that 
a subdivision
$$
\G
=
\G_1 \sqcup \ldots 
\sqcup \G_N
$$
is $\D$-fitted, if for any $j$
we have
$$\sum_{g_k\in\G_j} r_g=0.
$$
For any $\D$-fitted subdivision
$\ss$
consider the subgroup
$L(\ss)\sbs \Hom(\zz^n, \zz)$
consisting of all \ho s
$h:\zz^n\to \zz$ such that
$h~|~\G_i$ is constant for every $i$.
Then $L(\ss)$ is a full subgroup of $\Hom(\zz^n, \zz)$.
Observe that $L(\ss)$ is a proper subgroup since 
every $\G_i$ contains at least two elements.
It is clear that 
a \ho~
$p:\zz^n\to\zz^m$
belongs to $\Sigma$
if and only if $p$ is constant on each component
$\G_i$ 
of a $\D$-fitted subdivision of $\D$, that is,
$p\sqsubset L(\ss)$. $\qs$

\pa
\section{The twisted  Novikov Betti numbers}
\label{s:twi_nov_B}
\pa

Let $X$ be a finite connected CW complex, put $G=\pi_1(X)$. Let
 $R$ be an integral domain,
and $\eta:\pi_1(X)\to \GL(l,R)$ be a right representation 
(that is, $\eta$ is an antihomomorphism of groups).
Recall the definition of the twisted homological 
Betti numbers with coefficients in $\eta$:
$$
\b_k(X,\eta)
=
\dim_{\{R\}}
H_k\Big( \{ R\}^l \tens{\eta} C_*(\wi X) \Big).
$$
Starting with $\eta$ we can construct several other representations of $G$.
Let $n=\rk H_1(X,\zz)$ and denote by $\pi$ the projection 
$G\to H_1(G)/Tors\approx \zz^n$.
Let $L_n=R[\zz^n]$ denote the ring of Laurent polynomials in $n$ variables,
then the group $\zz^n$ can be identified with the group of 
units $L_n^{\times}\sbs \GL(1,L_n)$, and the \ho~ $\pi$ can be condsidered as 
a representation $G\to \GL(1,L_n)$. 
Denote by $\la\eta\ra$ the tensor product of $\pi$ and $\eta$, that is,
$$
\la\eta\ra(g)=\pi(g)\cdot \eta(g)\in \GL(l, L_n).
$$
Let $p:\zz^n\to\zz^m$ be a \ho.
Similarly to the above we can consider the tensor product of representations
$\eta$ and $p\circ\pi:G\to\GL(1, L_m)$.
This right representation will be denoted by $\la\eta\ra_p$.
Observe that
$\la\eta\ra_0=\eta, \ \ \la\eta\ra_{\Id}=\la\eta\ra. $
It is not difficult to see that for every $p:\zz^n\to\zz^m$ we have
\begin{equation}\label{f:different-eta}
\b_k(X,\eta) \geq \b_k(X,\la\eta\ra_p)
\geq
\b_k(X,\la\eta\ra).
\end{equation}
Theorem \ref{t:thins} of the previous section implies the following.
\bepr\lb{p:thins_twii}
Let $k\geq 0, q>0$.
Then there is a finite family of proper 
full subgroups $G_i\sbs \Hom(\zz^n, \zz)$
such that for $p\in\Hom(\zz^n, \zz^m)$
the condition
$$
\b_k(X,\la\eta\ra_p)
\geq 
\b_k(X,\la\eta\ra)
+q
$$
is equivalent to the following condition:
\
$p\sqsubset G_i$ \ for some $i$.
\enpr

Proceeding  to the Novikov homology, 
let us first recall the definition of the Novikov ring.
Let $H$ be a free abelian group; denote $\zz H$ by $\LL$.
Let $\mu:H\to\rr$ be a group \ho.
The Novikov completion $\wh\LL_\mu$ of the  ring $\LL$ 
\wrt~ $\mu$ is defined as the set of all series of the form $\l =\sum_g  n_g g$
(where $g\in H$  and  $n_g\in R$) satisfying the 
following  finiteness condition:
\begin{equation*}
\wh \LL_{\mu}
=\Big\{\l
~\Big|~ 
\forall
~C\in\rr,~ \
\text{the set }~\supp \l \cap \mu^{-1}\big([C, \infty [\big)
\text{ is finite  }\Big\}.
\end{equation*}
In general the ring 
$\wh \LL_{\mu}$
is rather complicated, 
however if $R=\zz$ and $\mu$ is a monomorphism, this ring is Euclidean 
by a theorem of J.-Cl. Sikorav (see \cite{PajitnovShrp}, Th. 1.4).
If $R$ is a field and  $\mu$ is a monomorphism, this ring is a field.

Let $\alpha:G\to\rr$ be a \ho. We can factor it as follows:

\bee\lb{f:factorr}
\xymatrix{G  \ar[r] & H_1(X,\zz)/Tors  \ar[r]^{ \hspace{1.3cm} \approx}  &
\zz^n  \ar[rr]  \ar[dr]_p &  & \rr \\
 & & & \zz^m \ar[ur]_{\wi\alpha}  & \\
}
\eee
where $p$ is an epimorphism and $\wi\alpha$ is a monomorphism.
We will denote the 
Novikov completion of the ring $R[\zz^m]$ \wrt~ $\wi\alpha$
 by $\wh L_{m, \alpha}$.
Denote by $\llc\eta\rrc$
the composition 
$$
G \arrr  {\eta_p} \GL(l, L_m) \arrinto \GL (l, \wh L_{m,\alpha});
$$
it is a right representation of $G$.
\bede\lb{d:twi_nov_hom}(\cite{GodaPajitnovTwiNov}, \cite{PajitnovTwiAlThurst})
The twisted homology of $X$ with respect to  the representation $\llc\eta\rrc$
is called  {\it the $\eta$-twisted 
Novikov homology of $X$ \wrt~ $\alpha$}.
The 
 $k$-th homological Betti number of $X$
\wrt~ $\llc\eta\rrc$  will be  denoted by $\wh b^\eta_k(X,\alpha)$
and called {\it the twisted Novikov Betti number of $X$ with respect to $\eta$}.
Thus we have 
%$$
%H_*\Big( \wh L_{m,\alpha}^l \tens{\eta_p}   C_*(\wi X)    \Big)
%$$
%is called {\it the $\eta$-twisted 
%Novikov homology of $X$ \wrt~ $\alpha$}.
%Put 
$$
\wh b^\eta_k(X,\alpha)
=
\dim_{\{\wh L_{m,\alpha}\}} H_k\Big( \{\wh L_{m,\alpha}\}^l \tens{\llc\eta\rrc} C_*(\wi X)     \Big).
$$
If $R=\zz$ so that 
$\wh L_{m,\alpha}$ is a principal ideal domain, 
we denote by $\wh q_k^\eta(X,\alpha)$ the torsion number of the module
$H_k\big(\wh L_{m,\alpha}^l \otimes C_*(\wi X)     \big)$.
\end{defi}

The geometric reasons to consider these completions 
of the module $C_*(\wi X)$ are as follows.
If $X$ is a compact manifold, and $\o$ is a closed 1-form on $M$ 
with non-degenerate zeroes,
denote by $\alpha$ the period homomorphism
$[\o]=\alpha : H_1(M,\zz)\to\rr$. 
The Morse-Novikov theory implies the following 
lower bounds on the number of zeroes of $\o$:
$$
m_k(\o)
\geq 
\frac 1l \Big( 
\wh b_k^\eta(X,\alpha)
+
\wh q_k^\eta(X,\alpha)
+
\wh q_{k-1}^\eta(X,\alpha)
\Big).
$$
(Here the base ring $R$ equals $\zz$.)
Compared to the other versions of the Novikov homology,
the twisted Novikov homology
has the advantage of being computable, and at the same time to
keep the information about the non-abelian structure of $G$
and related invariants. In a recent work 
\cite{FriedlVidussiVanTwiAl} S. Friedl and S. Vidussi
 proved that the twisted Novikov
homology detects fibredness of knots in $S^3$.
At present we will need only the simplest part 
of these  invariants, namely,
the twisted Novikov Betti numbers.
\bepr\lb{p:nov_localiz}
Let $\alpha:G\to \rr$ be a \ho~ and $\eta:G\to GL(l,R)$
a right representation. Then for every $k$ we have
$$
\b_k(X,\la\eta\ra_p) =\wh  b_k^\eta(X,\alpha)
$$
(where $p$ is obtained from the diagram 
\rrf{f:factorr}).
\enpr
\Prf 
The twisted Betti Novikov number in question
equals the dimension of the module
$$
H_k\big( L_m^l  \tens{\la\eta\ra_p}   C_*(\wi X)\big)\otimes \{\wh L_{m,\alpha}\}
$$
over the field of fractions $\{\wh L_{m,\alpha}\}$ of the Novikov ring.

The Betti number $b_k(X,\eta;p) $ is the dimension of the vector space
$$
H_k\big( L_m^l \tens{\la\eta\ra_p} C_*(\wi X)    \big)\otimes \{ L_{m}\}
$$
over the field of fractions $\{L_m\}$.
The inclusion $L_m \sbs\wh  L_{m,\alpha}$ extends to an inclusion of 
fields $\{L_m\} \sbs \{\wh L_{m,\alpha}\}$ and the result follows. $\qs$

\bede\lb{d:irr_deg}
For a \ho~ $\alpha : G\to\rr$ the subgroup $\alpha(G)$ is a free finitely generated
abelian group; its rank is called
{\it irrationality degree of $\alpha$}
and denoted $\Irr\alpha$. In particular $\Irr\alpha=1$ if and only if 
$\alpha$ is  a multiple of a \ho~ $G\to\zz$. We say that $\alpha$
is {\it maximally irrational } if $\Irr\alpha=\rk H_1(G)$.
\end{defi}
Observe that the irrationality degree of $\alpha$
equals the number $m$ from the diagram \rrf{f:factorr}.
If $\alpha$ is maximally irrational, then
the \ho~ 
$p$ in this diagram   is an
isomorphism. In this case the twisted Betti numbers 
$\wh b_k^\eta(X,\alpha)$
do not depend on $\alpha$.

\bede\lb{d:maxirr}
The number $\wh b_k^\eta(X,\alpha)$ where 
$\alpha$ is maximally irrational,
will be denoted by 
$\wh b_k^\eta(X)$. 
\end{defi}

The inequalities 
\rrf{f:different-eta} together with the proposition
\ref{p:nov_localiz} 
imply that
$$
\wh b^\eta_k(X,\alpha) \geq  \wh b_k^\eta(X).$$

Theorem \ref{t:thins} of the previous section implies the following.
\bepr\lb{p:thins_twi}
Let $k\geq 0, q>0$.
Then there is a finite family of proper 
full subgroups $G_i\sbs \Hom(\zz^n, \zz)$
such that 
%p\in\Hom(\zz^n, \zz^m)$
the condition
$$\wh b^\eta_k(X,\alpha)\geq \wh b^\eta_k(X)+q
$$
is equivalent to the following condition:
\
$\alpha\in \bigcup_i G_i\otimes\rr$. $\qs$
\enpr
\bere
\label{r:hom-cohomNov}
The results of this section have natural cohomology analogs.
Namely, given a representation $\r:G\to GL(l,R)$
one defines the representations
$\la\r\ra, \ \la\r\ra_p$, and the corresponding cohomological Betti numbers
$\b^k(X,\la\r\ra)$ and  $\b^k(X,\la\r\ra_p)$.  
One can define also the cohomological twisted Novikov numbers
$\wh b^k_\r(X,\alpha)$ and $\wh b^k_\r(X)$. The following lemma follows from 
Lemma \ref{l:hom_cohom}.
\bele\label{l:hom-cohom-nov}
Let $\r:G\to GL(l,R)$ be a representation;  put $\eta=\rho^*$.
We have
\begin{gather*}
\b^k\big(X,\la\r\ra\big)=\b_k\big(X,\la\eta\ra\big), \ \ \ \b^k\big(X,\la\r\ra_p\big)=\b_k\big(X,\la\eta\ra_p\big), \\
\wh b_k^\eta\big(X,\alpha\big)
=
\wh b^k_\r\big(X,\alpha\big)
, \ \ \ \wh b_k^\eta(X)=\wh b^k_\r(X).
 \hspace{2cm} \qs
\end{gather*}
\enle

\enre

\section{ Homology with local coefficients}
\lb{s:hom_loc}

Let us proceed now to the homology with local coefficients.
Let $\eta: G\to \GL(l,\cc)$ be a right representation of $G$
and $\a\in H^1(X,\cc)$.
The cohomology class $\a$ can be considered as a \ho~
$\a: G\to\cc$, which factors as follows 
\bee\lb{f:factor2}
\xymatrix{G  \ar[r] & H_1(X,\zz)/Tors  \ar[r]^{ \hspace{1.3cm} \approx}  &
\zz^n  \ar[rr]  \ar[dr]_p &  & \rr \\
 & & & \zz^m, \ar[ur]_{\wi\a}  & \\
}
\eee
where $p$ is an epimorphism and $\wi\a$ a monomorphism.
Here $m=\Irr\a$;
denote the coordinates of $p$ by $p_i :\zz^m\to\zz, \ 1\leq i \leq m$.
We have $\alpha=\sum_{i=1}^m\a_i p_i$,
and the numbers $\a_i\in\cc$ are linearly independent over $\qq$.
Recall  the exponential and formal exponential
deformations of $\eta$:
$$
\g_t:G\to \GL(l, \cc), \ \  
\g_t(g)=\eta(g)e^{t\langle\a, g\rangle}\in \cc \ \ \ ({\rm where}\ \  t\in\cc);
$$
$$
\wh\g(g)=\eta(g)e^{t\langle\a, g\rangle}\in \GL(l, \cc[[t]]).
$$
and the corresponding Betti numbers
$\b_k(X,\g_t), \ \b_k(X,\wh\g)$.
We will need a simple lemma.
\bele\lb{l:alg_indep}
Let $\a_1, \ldots , \a_m$ be complex numbers, linearly independent
over $\qq$. Then the power series $e^{\a_1 t}, \ldots , e^{\a_m t}\in \cc[[t]]$
are algebraically independent.
\enle
\Prf
If $P\in \cc[z_1, \ldots , z_m]$ is a polynomial such that
$P(e^{\a_1 t}, \ldots , e^{\a_m t})=0$, write
$P=\sum a_I t^I$
where the sum ranges over multiindices $I=(k_1, \ldots , k_m)\in \nn^m$.
Denote the string $(\a_1, \ldots , \a_m)$ by $\a$.
The  series $\zeta=P(e^{\a_1 t}, \ldots , e^{\a_m t})$ is then a finite
sum of exponential functions  of the form $a_Ie^{t\langle I, \a\rangle}$.
Observe that $\langle I, \a\rangle \not= \langle J, \a\rangle$
if $I\not= J$, since $\a_i$ are linearly independent over $\qq$.
Therefore $\z$ is a finite linear combination of exponential
functions  $e^{t\b_I}$
with pairwise different $\b_I$. Thus $\z=0$ implies $a_I=0$ for all $I$. $\qs$

\bepr\lb{p:equal}
There is a subset $S\sbs \cc$ consisiting of isolated points, such
that for every $t\in \cc\sm S$ 
and every $k$ we have

1) $\b_k(X,\g_t)=\b_k(X,\wh\g)= \b_k(X,\la\eta\ra_p)$

2) If $\a\in H^1(X,\rr)$
then 
$\b_k(X,\g_t)=\b_k(X,\wh\g) = \wh b_k^\eta(X,\a)$.
\enpr 
\Prf 
The right representation $\wh\g$ factors through $\zz^m$ as follows 
$$
G\arrr {p\circ \pi} \zz^m \arrr \G \Big(\cc[[t]]\Big)^*
$$
Let $(e_1, \ldots, e_m)$ denote the canonical basis in $\zz^m$;
then we have $\G(e_i)=e^{t\a_i}$.
The extension of $\G$ to a ring \ho~ 
$\zz[\zz^m]\to\cc[[t]]$ is injective by the previous lemma,
 and therefore can be further extended to a 
\ho~ 
$$
\ove \G : \{ L_m\} \to \cc((t))
$$
of the fraction fields. Therefore 

$$
\dim_{\{L_m\}}
H_k\Big( \{L_m\}^l \tens{\la\eta\ra_p} C_*(\wi X)
\Big)
=
\dim_{\cc((t))}
H_k\Big(
 \cc((t))^l
\tens{\G\circ\eta\circ\pi}
C_*(\wi X)
\Big)
$$
and the point 1) is proved. 
The point 2) follows from 1) and Proposition \ref{p:nov_localiz}. $\qs$

\bere\label{r:PapadimaSuciu}
 A particular case of this proposition
 corresponding to the vanishing Novikov Betti numbers 
was proved by S. Papadima and A. Suciu in 
\cite{PapadimaSuciuBNS}. 
\enre

The next proposition follows immediately.
\bepr\lb{p:thins_locC}
Let $k\geq 0, q>0$.
Then there is a finite family of proper 
full subgroups $G_i\sbs \Hom(\zz^n, \zz)$
such that 
%p\in\Hom(\zz^n, \zz^m)$
the condition
$$\b_k(X,\wh\g)\geq \b_k(X,\la\eta\ra)+q
$$
is equivalent to the following condition:
\
$\alpha\in \bigcup_i G_i\otimes\cc$.
\enpr

\section{Strongly formal spaces II : the jump loci}
\lb{s:form2}

Let $X$ be a manifold, $\r:\pi_1(X)\to \GL(n,\cc)$
a representation, and $\a\in H^1(X,\cc)$.
In this section we assume that the spectral sequences
$\EEEE^*_r$ and $\WWWW^*_r$ associated to 
the triple $(X,\r,\a)$ 
degenerate at their second term, so that
$$
E^*_r(X,\r,\a)
\approx 
W^*_r(X,\r,\a)
\approx
\HH^*\Big( H^*(X,\r), \a\Big).
$$
According to Section \ref{s:form_hypform}
this holds, in particular, when $X$ is a K\"ahler manifold
and $\r$ a semi-simple representation.
%The following theorem follows from Theorem \ref{t:compar_manif},
%together with Propositions \ref{p:nov_localiz}, \ref{p:equal}.
Recall from Section \ref{s:spec_twis}
the exponential deformation $\g_t$ and the formal exponential deformation $\wh\g$ of the
representation $\r$.

\beth\lb{t:form2}
\been\item
$\b^k(X,\wh\g) = \b^k(X,\g_{gen}) 
= \BB^k(H^*(X,\r), \a).$
\item If $\a\in H^1(X,\rr)$ then 
$$\wh b^k_\r(X,\a) = \b^k(X,\wh\g) 
= \BB^k(H^*(X,\r), \a). \qquad $$
\enen
\enth
\Prf
The first point is a consequence of the degeneracy of the two spectral sequences above.
To prove the second point, consider the right representation $\eta=\r^*$,
and the corresponding formal exponential deformation $\wh\g^*$ of $\eta$.
We have 
$
\wh b^k_\r(X,\a)
=
\wh b_k^\eta(X,\a)
$
by Lemma \ref{l:hom-cohom-nov}.
Further, 
$
\wh b_k^\eta(X,\a)
=
\b_k(X, \wh\g^*)$
by Proposition \ref{p:equal}.
Finally 
$\b_k(X, \wh\g^*)
=
\b^k(X, \wh\g)$
by Corollary \ref{c:hom_cohom}.
The proof of the theorem is over. $\qs$

The proof of the following proposition,
concerning the jump
loci of the Betti numbers,
is done on similar lines, using 
Proposition \ref{p:thins_locC}
and Proposition \ref{p:thins_twii}.

\bepr\lb{t:form2jump}
Let $k\geq 0, q>0$.
Then there is a finite family of proper 
full subgroups $G_i\sbs \Hom(\zz^n, \zz)$
such that each of the following conditions \rrf{f:cond1}, \rrf{f:cond2}
\begin{equation}\lb{f:cond1}
\b^k(X,\wh\g) \geq \b^k(X,\la\r\ra) + q 
\end{equation}
\begin{equation}\lb{f:cond2}
\BB^k(H^*(X,\r), \a) \geq \b^k(X,\la\r\ra) + q
\end{equation}
is equivalent to the condition
$\a\in \bigcup_i G_i\otimes\cc$.

For 
$\a\in H^1(X,\rr)$
each of the following conditions \rrf{f:cond3}, \rrf{f:cond4}
\begin{equation}\lb{f:cond3}
\wh b_\r^k(X,\alpha)\geq \b^k(X,\la\r\ra)+q; 
\end{equation}
\begin{equation}\lb{f:cond4}
\BB^k(H^*(X,\r), \a) \geq \b^k(X,\la\r\ra) + q
\end{equation}
is equivalent to the condition
$\a\in \bigcup_i G_i\otimes\rr$. $\qs$
\enpr

\section{Acknowledgements}

This work was supported by World Premier International Research Center
Initiative (WPI Program), MEXT, Japan and by
JSPS Grants-in-Aid for Scientific Research 23340014.
A part of this work was accomplished
while the second author was staying at Graduate School of
Mathematical Sciences, the University of Tokyo, in 2011.
The second author gratefully  acknowledges 
the support by the Program for Leading Graduate
Schools, MEXT, Japan, and thanks the Graduate School of
Mathematical Sciences for warm hospitality.

The authors thank Hisashi Kasuya for several valuable discussions
concerning his work and the subject of the present paper.

During the preparation of the paper, the authors became 
aware of the e-print  
\cite{DimcaPapadimaNonAbJL}
by A. Dimca and S. Papadima,
on related subjects. Theorem \ref{t:formal_top}
of the present paper is also proved by a completely 
different method in \cite{DimcaPapadimaNonAbJL},
Theorem E.

The second author thanks Fedor Bogomolov for many valuable discussions
and for sharing his insight about several questions of algebraic geometry.

The authors thank the anonymous referees for the helpful comments
and suggestions which lead to  the improvement of the manuscript.

\input ref.tex
\end{document}

%% file: auxxx.tex
\renewcommand{\a}{\alpha}
\renewcommand{\b}{\beta}
\newcommand{\g}{\gamma}

\newcommand{\z}{\zeta}
\renewcommand{\t}{\theta}
\renewcommand{\l}{\lambda}

\renewcommand{\r}{\rho}

\renewcommand{\o}{\omega}

\newcommand{\G}{\Gamma}
\newcommand{\D}{\Delta}

\renewcommand{\L}{\Lambda}

\renewcommand{\AA}{{\mathcal A}}
\newcommand{\BB}{{\mathcal B}}

\newcommand{\HH}{{\mathcal H}}

\newcommand{\KK}{{\mathcal K}}
\newcommand{\LL}{{\mathcal L}}
\newcommand{\MM}{{\mathcal M}}
\newcommand{\NN}{{\mathcal N}}

\newcommand{\QQ}{{\mathcal Q}}
\newcommand{\RR}{{\mathcal R}}
\renewcommand{\SS}{{\mathcal S}}
\newcommand{\TT}{{\mathcal T}}

\newcommand{\ZZ}{{\mathcal Z}}

\newcommand{\cc}{{\mathbb{C}}}

\newcommand{\kk}{{\mathbb{K}}}

\newcommand{\nn}{{\mathbb{N}}}

\newcommand{\qq}{{\mathbb{Q}}}
\newcommand{\rr}{{\mathbb{R}}}
\renewcommand{\ss}{{\mathbb{S}}}

\newcommand{\zz}{{\mathbb{Z}}}

\newcommand{\RRR}{{\mathbf{R}}}

\newcommand{\CCCC}{{\mathscr{C}} }

\newcommand{\EEEE}{{\mathscr{E}} }

\newcommand{\WWWW}{{\mathscr{W}}}

\newcommand{\llc}{\langle\langle}
\newcommand{\rrc}{\rangle\rangle}

\newcommand{\la}{\langle}
\newcommand{\ra}{\rangle}

\newcommand{\Ker}{\text{\rm Ker }}

\newcommand{\Hom}{\text{\rm Hom}}

\newcommand{\rk}{\text{\rm rk}}
\renewcommand{\Im}{\text{\rm Im }}
\newcommand{\supp}{\text{\rm supp }}

\newcommand{\Id}{\text{\rm Id}}

\newcommand{\Irr}{\text{\rm Irr}}

\newcommand{\GL}{\text{\rm GL}}

\newcommand{\bere}{\begin{rema}}
\newcommand{\bede}{\begin{defi}}

\renewcommand{\beth}{\begin{theo}}
\newcommand{\bele}{\begin{lemm}}
\newcommand{\bepr}{\begin{prop}}
\newcommand{\beeq}{\begin{equation}}
\newcommand{\bega}{\begin{gather}}
\newcommand{\begaa}{\begin{gather*}}
\newcommand{\been}{\begin{enumerate}}

\newcommand{\bedee}{\begin{defii}}
\newcommand{\bethh}{\begin{theoo}}
\newcommand{\belee}{\begin{lemmm}}
\newcommand{\beprr}{\begin{propp}}

\newcommand{\beco}{\begin{coro}}

\newcommand{\beal}{\begin{aligned}}

\newcommand{\enre}{\end{rema}}

\newcommand{\enco}{\end{coro}}
\newcommand{\enpr}{\end{prop}}
\newcommand{\enth}{\end{theo}}
\newcommand{\enle}{\end{lemm}}
\newcommand{\enen}{\end{enumerate}}
\newcommand{\enga}{\end{gather}}
\newcommand{\engaa}{\end{gather*}}
\newcommand{\eneq}{\end{equation}}
\newcommand{\enal}{\end{aligned}}

\newcommand{\bq}{\begin{equation}}
\newcommand{\bqq}{\begin{equation*}}

\renewcommand{\leq}{\leqslant}
\renewcommand{\geq}{\geqslant}

\newcommand{\wi}{\widetilde}

\newcommand{\ove}{\overline}

\newcommand{\wh}{\widehat}

\newcommand{\sps}{\supset}

\newcommand{\sm}{\setminus}

\newcommand{\sbs}{\subset}

\newcommand{\tens}[1]{\underset{#1}{\otimes}}

\newcommand{\wrt}{with respect to}
\newcommand{\ho}{homomorphism}

\newcommand{\Prf}{{\it Proof.\quad}}

\newcommand{\smo}{C^{\infty}}

\newcommand{\chart}{\Phi_p:U_p\to B^n(0,r_p)}
\newcommand{\atlas}{\{\Phi_p:U_p\to B^n(0,r_p)\}_{p\in S(f)}}

\newcommand{\pr}{\partial}

\newcommand{\qs}{\hfill\square}

\newcommand{\pa}{\vskip0.1in}